# STABILITY OF NONLINEAR FILTERS IN NONMIXING CASE

By Pavel Chigansky and Robert Liptser

*University of Tel Aviv*

The nonlinear filtering equation is said to be stable if it "forgets" the initial condition. It is known that the filter might be unstable even if the signal is an ergodic Markov chain. In general, the filtering stability requires stronger signal ergodicity provided by the, so called, mixing condition. The latter is formulated in terms of the transition probability density of the signal. The most restrictive requirement of the mixing condition is the uniform positiveness of this density. We show that it might be relaxed regardless of an observation process structure.

**1. Introduction and the main result.** This paper addresses the stability problem of the nonlinear filtering equation with respect to its initial condition. We consider a homogeneous ergodic Markov chain $(X_n)_{n\geq 0}$ with values in $\mathbb{S}\subseteq\mathbb{R}^d$ regarded as a signal to be filtered from observation of $(Y_n)_{n\geq 1}$, $Y_n\in\mathbb{R}^p$.

Denote $M(dx)$, $K(x;dy)$ and $N(dx)$ the invariant measure, the transition probability kernel and the distribution of $X_0$, respectively, and describe the dependence of observation and signal processes via the conditional distribution $\Gamma(x;dy)$: for any Borel measurable set $A$ from $\mathbb{R}^p$ and $X_{[0,n]}$, $Y_{[1,n]}$ the $\sigma$-algebras generated by $\{X_0,\dots,X_n\}$, $\{Y_1,\dots,Y_n\}$

$$P(Y_n\in A|X_{[0,n]}\vee Y_{[1,n-1]})=P(Y_n\in A|X_n)=\int_A\Gamma(X_n,dy)\qquad\text{a.s.}\tag{1.1}$$

An example of such type dependence is

$$Y_n=h(X_n,\xi_n),$$

with a bounded measurable function $h$ and i.i.d. sequence of random variables $(\xi_n)_{n\geq 1}$ independent of $(X_n)_{n\geq 0}$.

---

We assume that there exist reference measures $\psi(dy)$ on $\mathbb{S}$ and $\theta(dy)$ on $\mathbb{R}^p$ such that

$$K(x;dy)=\lambda(x,y)\psi(dy),$$
$$N(dx)=\nu(x)\psi(dx),$$
$$\Gamma(x;dy)=\gamma(x,y)\theta(dy),$$

where the densities $\lambda(x,y)$ and $\gamma(x,y)$ are measurable nonnegative functions. Then, obviously, $M(dx)=m(x)\psi(dx)$ and $P(X_n\in dx|Y_{[1,n]})$ has (a.s.) a density $\pi_n^\nu(x)$ with respect to $\psi$: for any Borel set $A$ from $\mathbb{S}$,

$$P(X_n\in A|Y_{[1,n]})=\int_A \pi_n^\nu(x)\psi(dx).$$

It is well known that $\pi_n(x)$, $n\ge 1$ solves the, so called, filtering equation (recursive Bayes' formula)

$$\pi_n^\nu(x)=\frac{\gamma(x,Y_n)\int_{\mathbb{S}}\lambda(z,x)\pi_{n-1}^\nu(z)\psi(dz)}{\int_{\mathbb{S}}\int_{\mathbb{S}}\gamma(u,Y_n)\lambda(v,u)\pi_{n-1}^\nu(v)\psi(dv)\psi(du)} \tag{1.2}$$

subject to $\pi_0^\nu(x)=\nu(x)$.

If the density $\nu(x)$ is not completely known, in practice the filtering equation is initialized by another density $\beta(x)$ corresponding to some probability distribution $B(dx)$ on $\mathbb{S}$ with respect to $\psi$. Though this substitution seems natural, the fraction in the right-hand side of (1.2) may not be well defined. Nevertheless, if $N\ll B$ [henceforth, $N\ll B$ is assumed and $\frac{dN}{dB}(x)=:\frac{\nu}{\beta}(x)$], the recurrent equation (1.2) subject to $\beta$ makes sense $\psi$-a.s.

If $\nu$ is replaced by $\beta$, we get "wrong" filtering density $\pi_n^{\beta\nu}(x)$ defined as follows: $\pi_0^{\beta\nu}(x)=\beta(x)$ and

$$\pi_n^{\beta\nu}(x)=\frac{\gamma(x,Y_n)\int_{\mathbb{S}}\lambda(z,x)\pi_{n-1}^{\beta\nu}(z)\psi(dz)}{\int_{\mathbb{S}}\int_{\mathbb{S}}\gamma(u,Y_n)\lambda(v,u)\pi_{n-1}^{\beta\nu}(v)\psi(dv)\psi(du)}.$$

Since $\pi_n^{\beta\nu}$ differs from $\pi_n^\nu$, the use of $\pi_n^{\beta\nu}\langle f\rangle:=\int_{\mathbb{S}}f(x)\pi_n^{\beta\nu}(x)\psi(dx)$ as an estimate for $f(X_n)$, with bounded $f$ or $Ef^2(X_n)<\infty$, does not guarantee the mean square optimality, that is,

$$E\big(f(X_n)-\pi_n^\nu\langle f\rangle\big)^2<E\big(f(X_n)-\pi_n^{\beta\nu}\langle f\rangle\big)^2$$

even as $n\nearrow\infty$.

The filtering equation, or filter, is said to be stable, if for any bounded measurable function $f$,

$$\lim_{n\to\infty}E|\pi_n^\nu\langle f\rangle-\pi_n^{\beta\nu}\langle f\rangle|=0, \tag{1.3}$$

so that, from the practical point of view, the use of $\pi_n^{\beta\nu}$ for stable filters makes sense.

Assume that $\beta$ is chosen and consider a new pair of processes $(X^\beta_n, Y^\beta_n)$ characterized by the same kernels $K$ and $\Gamma$ and $\beta$ as the distribution density of $X^\beta_0$. It can be shown (similar to Proposition 2.1, [2]) that (1.3) holds provided that

$$E\Big(\frac{\nu}{\beta}(X^\beta_0)\Big|Y^\beta_{[1,\infty)}\Big) = \lim_{n\to\infty} E\Big(\frac{\nu}{\beta}(X^\beta_0)\Big|X^\beta_n, Y^\beta_{[1,n]}\Big) \qquad \text{a.s.} \tag{1.4}$$

So, (1.3) follows from $Y^\beta_{[1,\infty)} = \bigcap_{n\ge 0}\{Y^\beta_{[1,\infty)} \vee X^\beta_{[n,\infty)}\}$, since by the Markov property and the martingales convergence theorem,

$$\lim_{n\to\infty} E\Big(\frac{\nu}{\beta}(X^\beta_0)\Big|X_n, Y^\beta_{[1,n]}\Big) = E\Big(\frac{\nu}{\beta}(X^\beta_0)\Big|\bigcap_{n\ge 0}\{Y^\beta_{[1,\infty)} \vee X^\beta_{[n,\infty)}\}\Big).$$

Unfortunately, in general, $Y^\beta_{[1,\infty)} \subsetneq \bigcap_{n\ge 0}\{Y^\beta_{[1,\infty)} \vee X^\beta_{[n,\infty)}\}$ even though the ergodic property of the Markov chain implies $\bigcap_{n\ge 0} X^\beta_{[n,\infty)} = (\varnothing, \Omega)$ a.s.

In several publications $Y^\beta_{[1,\infty)} = \bigcap_{n\ge 0}\{Y^\beta_{[1,\infty)} \vee X^\beta_{[n,\infty)}\}$ is implicitly declared and then (1.3) is derived (for more details, see [2] and citations therein). Somewhat contradicting the intuition, the filter does not automatically inherit (1.3) from the ergodic property of signal. Corresponding counterexamples can be found in [2, 7, 8] (see also Example 5.1).

The validity of (1.3) sometimes crucially depends on the structure of $\gamma(x, y)$. For instance, in order to have (1.3), $\gamma(x, y)$ is required to be compactly supported (see [4]), or to have tails with certain decay rate depending on the signal model (see [1, 5, 9, 10, 16]).

In view of the aforementioned facts, it is interesting to determine the conditions providing (1.3) regardless of $\gamma(x, y)$. One of them is the *mixing condition* introduced by Atar and Zeitouni [1] and Del Moral and Guionnet [6]:

$$0 < \lambda_* \le \lambda(x, y) \le \lambda^* \tag{1.5}$$

which provides

$$\limsup_{n\to\infty} \frac{1}{n} \log \|\pi^\nu_n - \pi^{\beta\nu}_n\| \le -\frac{\lambda_*}{\lambda^*} \qquad \text{a.s.,} \tag{1.6}$$

where

$$\|\pi^{\beta\nu}_n - \pi^\nu_n\| = \int_{\mathbb{S}} |\pi^{\beta\nu}_n(x) - \pi^\nu_n(x)|\psi(dx)$$

is the total variation norm. Clearly (1.6) implies (1.3).

If after all the condition $\lambda^* < \infty$ is reasonable, the condition $\lambda_* > 0$ is quite restrictive in many applications.

In this paper we replace $\lambda_*$ by

$$\lambda_\diamond = \int_{\mathbb{S}} \operatorname*{ess\,inf}_{y\in\mathbb{S}} \lambda(x,y)m(x)\psi(dx) \tag{1.7}$$

(hereafter, ess inf and ess sup are taken with respect to $\psi$).

THEOREM 1.1. *Assume $0 < \lambda_\diamond$ and $\lambda^* < \infty$. Then, for any $\nu(x)$ and any $\beta(x) \ge \beta_* > 0$, the relation (1.4) holds and*

$$\limsup_n \frac{1}{n} \log \|\pi_n^{\beta\nu} - \pi_n^\nu\| \le -\frac{\lambda_\diamond}{\lambda^*} \qquad a.s. \tag{1.8}$$

REMARK 1.1. In contrast to $\lambda(x,y) \ge \lambda_* > 0$, the assumptions of Theorem 1.1 admit $\lambda(x,y)$ to be zero on some region of the state space.

EXAMPLE 1.1. Let $X_n$ be the Markov chain with values in $\mathbb{S} = \{1,\dots,d\}$ and $\Lambda$ be its transition probability matrix with entries

$$\lambda_{ij} = P(X_n = j | X_{n-1} = i).$$

Here $\lambda_{ij}$ is the transition probability density with respect to the Dirac measure $\psi$ supported at points $j = 1,\dots,d$.

Assume $\lambda_* = \min_{ij} \lambda_{ij} = 0$, but there is an integer $r > 1$ such that all entries of $\Lambda^r$ are positive. Then the Markov chain is ergodic (see, e.g., [14]) and its invariant measure possesses density $m$ with respect to $\psi$ with positive atoms $m(i)$, $i = 1,\dots,d$. Then

$$\lambda_\diamond = \sum_{i=1}^d \min_j \lambda_{ij} m(i) > 0,$$

if at least one row of $\Lambda$ has positive entries.

REMARK 1.2. The relation $m(y) = \int_{\mathbb{S}} \lambda(x,y)m(x)\psi(dx)$ implies

$$\lambda_\diamond \le m(x) \le \lambda^*. \tag{1.9}$$

REMARK 1.3. The requirement $\beta \ge \beta_* > 0$ is not restrictive, since the choice of $\beta$ is flexible.

Notice also that $\beta \ge \beta_*$ implies $N \ll B$ with

$$\frac{dN}{dB}(x) = \frac{\nu(x)}{\beta(x)}.$$

Our proof of (1.8) differs from these in [1] and [6]. Before getting into details, we describe the main idea of the proof, where the *backward* conditional density

$$\rho_n(u,x) = \frac{P(X_0 \in du | X_n = x, Y_{[1,n]})}{\psi(du)}$$

plays a crucial role.

We show in Section 4.3 that $|\rho_n(u,x') - \rho_n(u,x'')| \to 0$, $n \to \infty$, for any $x'$, $x''$ which allows the claim that

$$\lim_{n\to\infty} E\Big(\frac{\nu}{\beta}(X_0^\beta)|X_n^\beta \vee Y_{[1,n]}^\beta\Big)$$

is a $Y_{[1,\infty)}$-measurable random variable and in turn the validity of (1.4).

The proof of the second statement of Theorem 1.1 uses an upper bound for the rate, in $n$, of $|\rho_n(u,x') - \rho_n(u,x'')| \to 0$. When $\lambda_* > 0$, this rate is $\lambda_*/\lambda^*$, while for $\lambda_\diamond > 0$ the rate is random (Lemma 3.2) and controlled by $\pi_n(x)$. The upper bound in (1.8) follows from the law of large numbers for the conditional expectations (Theorem 2.2), which is derived with the help of geometric ergodicity (Theorem 2.1) being valid under the assumptions of Theorem 1.1.

## 2. Geometric ergodicity and law of large numbers for conditional expectations.

2.1. *Geometric ergodicity.* Let $\lambda^{(n)}(u,x)$ be the $n$-steps transition probability density with respect to $\psi$.

The Markov chain is geometrically ergodic, if there exist constants $C > 0$ and $0 < r < 1$ such that for any $u$,

$$\int_{\mathbb{S}} |\lambda^{(n)}(u,x) - m(x)|\psi(dx) \le Cr^n, \qquad n \ge 1.$$

THEOREM 2.1. *Assume $0 < \lambda_\diamond$ and $\lambda^* < \infty$. Then the geometric ergodicity holds with*

$$C = \frac{(\lambda^*)^2}{\lambda_\diamond(\lambda^* - \lambda_\diamond)} \quad \textit{and} \quad r = 1 - \frac{\lambda_\diamond}{\lambda^*}.$$

PROOF. Consider the stationary Markov chain $X_n$, that is, when $X_0$ has distribution density $m(x)$ with respect to $\psi$.

By the Bayes formula, one verifies that $P(X_0 \in du | X_n = x)$ has a density $q_n(u,x)$ with respect to $\psi$

$$q_1(u,x) = \frac{\lambda(u,x)m(u)}{m(x)} \tag{2.1}$$

and

$$q_n(u,x)=\frac{\int_{\mathbb{S}}\lambda(x',x)q_{n-1}(u,x')m(x')\psi(dx')}{m(x)}. \tag{2.2}$$

Let us show now that

$$\begin{aligned} &m(u)[\lambda^{(n)}(u,x)-m(x)]\\ &\qquad=\int_{\mathbb{S}}[q_n(u,x)-q_n(u,x')]m(x)m(x')\psi(dx'). \end{aligned} \tag{2.3}$$

For arbitrary bounded and measurable functions $g$ and $f$, write

$$\begin{aligned} &Eg(X_0)f(X_n)-Eg(X_0)Ef(X_n)\\ &\qquad=\int_{\mathbb{S}}\int_{\mathbb{S}}g(u)f(x)[\lambda^{(n)}(u,x)-m(x)]m(u)\psi(du)\psi(dx) \end{aligned} \tag{2.4}$$

and

$$\begin{aligned} Eg(X_0)f(X_n)&=\int_{\mathbb{S}}\int_{\mathbb{S}}g(u)f(x)q_n(u,x)m(x)\psi(du)\psi(dx)\\ &=\int_{\mathbb{S}}\int_{\mathbb{S}}\int_{\mathbb{S}}g(u)f(x)q_n(u,x)m(x)m(x')\psi(du)\psi(dx)\psi(dx'). \end{aligned} \tag{2.5}$$

Obviously, (2.5) provides

$$\begin{aligned} &Eg(X_0)f(X_n)-Eg(X_0)Ef(X_n)\\ &\qquad=\int_{\mathbb{S}}\int_{\mathbb{S}}\int_{\mathbb{S}}g(u)f(x)[q_n(u,x)-q_n(u,x')]\\ &\qquad\qquad\times m(x)m(x')\psi(du)\psi(dx)\psi(dx'). \end{aligned} \tag{2.6}$$

So, by arbitrariness of $g$ and $f$, (2.6) and (2.4) imply (2.3).

Owing to (1.9), the relation (2.3) provides

$$\begin{aligned} &|\lambda^{(n)}(u,x)-m(x)|\\ &\qquad=\frac{1}{m(u)}\int_{\mathbb{S}}|q_n(u,x)-q_n(u,x')|m(x)m(x')\psi(dx'). \end{aligned} \tag{2.7}$$

Introduce $\overline{q}_n(u)=\operatorname{ess\,sup}_{r\in\mathbb{S}}q_n(u,r)$, $\underline{q}_n(u)=\operatorname{ess\,inf}_{r\in\mathbb{S}}q_n(u,r)$ and set

$$\triangle_n(u)=\overline{q}_n(u)-\underline{q}_n(u).$$

By (2.7),

$$\int_{\mathbb{S}}|\lambda^{(n)}(u,x)-m(x)|\psi(dx)\le\frac{1}{m(u)}\triangle_n(u). \tag{2.8}$$

We show now that

$$\triangle_n(u) \le m(u)\frac{\lambda^*}{\lambda_\diamond}\Big(1-\frac{\lambda_\diamond}{\lambda^*}\Big)^{n-1}, \qquad n \ge 1. \tag{2.9}$$

For $n=1$, this estimate is obvious, whereas $q_1(u,x) \le m(u)\lambda^*/\lambda_\diamond$, that is, $\overline{q}_1(u) \le m(u)\lambda^*/\lambda_\diamond$, and $\triangle_1(u) \le \overline{q}_1(u)$. In order to establish this estimate for $n \ge 2$, we derive at first a recurrent inequality

$$\triangle_n(u) \le \Big(1-\frac{\lambda_\diamond}{\lambda^*}\Big)\triangle_{n-1}(u). \tag{2.10}$$

For any $v', v''$, we have

$$\overline{q}_{n-1}(u) = \frac{\int_{\mathbb{S}} m(x')\lambda(x',v')\overline{q}_{n-1}(u)\psi(dx')}{m(v')},$$

$$\underline{q}_{n-1}(u) = \frac{\int_{\mathbb{S}} m(x')\lambda(x',v'')\underline{q}_{n-1}(u)\psi(dx')}{m(v'')}$$

and, by (2.2),

$$\begin{aligned} &q_n(u,v') - q_n(u,v'') \\ &\quad = \frac{\int_{\mathbb{S}} \lambda(x',v')q_{n-1}(u,x')m(x')\psi(dx')}{m(v')} \\ &\qquad - \frac{\int_{\mathbb{S}} \lambda(x',v'')q_{n-1}(u,x')m(x')\psi(dx')}{m(v'')}. \end{aligned}$$

Then

$$\begin{aligned} &q_n(u,v') - q_n(u,v'') \\ &\quad = \overline{q}_{n-1}(u) - \underline{q}_{n-1}(u) \\ &\qquad - \frac{\int_{\mathbb{S}} m(x')\lambda(x',v')[\overline{q}_{n-1}(u) - q_{n-1}(u,x')]\psi(dx')}{m(v')} \\ &\qquad - \frac{\int_{\mathbb{S}} m(x')\lambda(x',v'')[q_{n-1}(u,x') - \underline{q}_{n-1}(u)]\psi(dx')}{m(v'')}. \end{aligned}$$

This equality, with ($0/0 = 0.5$ is understood here)

$$\alpha_n(u,x') := \frac{q_n(u,x') - \underline{q}_n(u)}{\triangle_n(u)} \quad \text{and} \quad 1-\alpha_n(u,x') := \frac{\overline{q}_n(u) - q_n(u,x')}{\triangle_n(u)},$$

is transformed to

$$\begin{aligned} &q_n(u,v') - q_n(u,v'') \\ &\quad = \triangle_{n-1}(u)\Big(1 - \int_{\mathbb{S}} m(x')\Big\{\frac{\lambda(x',v')}{m(v')}[1-\alpha_{n-1}(u,x')] \end{aligned}$$

$$+\frac{\lambda(x',v'')}{m(v'')}\alpha_{n-1}(u,x')\Big\}\psi(dx')\Big).$$

Owing to $0\le\alpha_n(u,x')\le 1$ and $m(x)\le\lambda^*$, we get

$$\begin{aligned}&\int_{\mathbb{S}} m(x')\Big\{\frac{\lambda(x',v')}{m(v')}[1-\alpha_{n-1}(u,x')]+\frac{\lambda(x',v'')}{m(v'')}\alpha_{n-1}(u,x')\Big\}\psi(dx')\\ &\quad\ge\int_{\mathbb{S}}\frac{m(x')}{\lambda^*}[\lambda(x',v')\wedge\lambda(x',v'')]\psi(dx')\\ &\quad\ge\frac{1}{\lambda^*}\int_{\mathbb{S}} m(x')\operatorname*{ess\,inf}_r\lambda(x',r)\psi(dx')=\frac{\lambda_\diamond}{\lambda^*}.\end{aligned}$$

Thus, $q_n(u,v')-q_n(u,v'')\le(1-\lambda_\diamond/\lambda^*)\triangle_{n-1}(u)$, and (2.10) holds true by arbitrariness of $v'$ and $v''$.

The iteration of (2.10) and the estimate $\triangle_1(u)\le m(u)\lambda^*/\lambda_\diamond$ give (2.9). So, by (2.8) the result holds. □

2.2. *The law of large numbers.* The result of this section is valid for any distribution of $X_0$, obeying density with respect to $\psi$.

THEOREM 2.2. *Assume $\lambda_\diamond>0$ and $\lambda^*<\infty$. Then, for any bounded measurable function $f$,*

$$\lim_{n\to\infty}\frac{1}{n}\sum_{k=1}^n\pi_{k-1}\langle f\rangle=\int_{\mathbb{S}} f(x)m(x)\psi(dx)\qquad a.s.$$

PROOF. Set $f^\circ(x)=f(x)-\int_{\mathbb{S}} f(z)m(z)\psi(dz)$ and notice that the desired statement is valid if

$$\lim_{n\to\infty}\frac{1}{n}\sum_{k=1}^n\pi_{k-1}\langle f^\circ\rangle=0\qquad\text{a.s.}\tag{2.11}$$

Let us introduce the Poisson equation

$$g(x)=f^\circ(x)+\int_{\mathbb{S}} g(y)\lambda(x,y)\psi(dy),$$

whose solution

$$g(x)=f^\circ(x)+\sum_{n=1}^\infty\int_{\mathbb{S}} f^\circ(y)\lambda^{(n)}(x,y)\psi(dy)$$

is well defined and bounded by virtue of the geometric ergodicity (see Theorem 2.1).

Set $\xi_n = \pi_n\langle g\rangle - \pi_{n-1}\langle g\rangle + \pi_{n-1}\langle f^\circ\rangle$ and notice that $\xi_n$'s forms a bounded martingale difference. The boundedness is obvious; $E(\xi_n|Y_{[1,n-1]}) = 0$, a.s., $n \ge 1$, is valid by the following Poisson equation:

$$\begin{aligned} E(\xi_n|Y_{[1,n-1]}) &= E(g(X_n) - g(X_{n-1}) + f^\circ(X_{n-1})|Y_{[1,n-1]}) \\ &= E\biggl(\int_{\mathbb{S}} g(z)\lambda(X_{n-1},z)\psi(dz) - g(X_{n-1}) + f^\circ(X_{n-1})\Big|Y_{[1,n-1]}\biggr) \\ &= 0. \end{aligned}$$

So, with the martingale $M_n = \sum_{k=1}^n \xi_k$, we have

$$\pi_n\langle g\rangle = \pi_0\langle g\rangle - \sum_{k=1}^n \pi_{k-1}\langle f^\circ\rangle + M_n.$$

Consequently,

$$\frac{1}{n}\sum_{k=1}^n \pi_{k-1}\langle f^\circ\rangle = \frac{1}{n}M_n - \frac{1}{n}(\pi_n\langle g\rangle - \pi_0\langle g\rangle).$$

With bounded $g$, $\pi_n\langle g\rangle$ can be chosen bounded too. Hence,

$$n^{-1}(\pi_n\langle g\rangle - \pi_0\langle g\rangle) \underset{n\to\infty}{\to} 0.$$

Also $M_n/n \to 0$, a.s. $n \to \infty$, since the increments of martingale $M_n$ are bounded (see Theorem 4, Chapter VII, Section 5 in [17]).

Thus, (2.11) holds. $\square$

## 3. The backward distribution $P(X_0 \in du|X_n \equiv x, Y_{[1,n]})$.

The results of this section are valid for any distribution of $X_0$ obeying density with respect to $\psi$, say $\vartheta$, bounded below by a constant $\vartheta_* > 0$.

3.1. *Recurrent equation for the density.* We show that the backward distribution possesses density $\rho_n(u,x)$ with respect to $\psi$.

Lemma 3.1. *Under the assumptions of Theorem 1.1 $\rho_n(u,x)$ exists and is defined as*

$$\begin{aligned} \rho_1(u;x) &= \frac{\lambda(u;x)\vartheta(u)}{\int_{\mathbb{S}}\lambda(v;x)\vartheta(v)\psi(dv)}, \\ \rho_n(u;x) &= \frac{\int_{\mathbb{S}}\lambda(x';x)\rho_{n-1}(u;x')\pi_{n-1}(x')\psi(dx')}{\int_{\mathbb{S}}\lambda(x';x)\pi_{n-1}(x')\psi(dx')}, \qquad n \ge 2. \end{aligned} \tag{3.1}$$

Proof. For $n = 1$, the formula is obvious.

By (1.1), $E(h(X_0)|X_n, Y_{[1,n]}) = E(h(X_0)|X_n, Y_{[1,n-1]})$ for any bounded measurable $h$ and $n \ge 2$. We apply the induction method. With bounded and measurable functions $g(x), f(x),\ x\in\mathbb{S}$, and $H(y_1,\dots,y_{n-1}),\ y_i\in\mathbb{R}^p$, write

$$\begin{aligned}
&Eg(X_0)H(Y_1,\dots,Y_{n-1})f(X_n)\\
&\qquad = Eg(X_0)H(Y_1,\dots,Y_{n-1})\int_{\mathbb{S}} f(x)\lambda(X_{n-1},x)\psi(dx)\\
&\qquad = EH(Y_1,\dots,Y_{n-1})\int_{\mathbb{S}}\int_{\mathbb{S}}\int_{\mathbb{S}} g(u)f(x)\lambda(x',x)\rho_{n-1}(u;x')\\
&\qquad\qquad\times \pi_{n-1}(x')\psi(dx)\psi(dx')\psi(du).
\end{aligned}\tag{3.2}$$

On the other hand,

$$\begin{aligned}
&Eg(X_0)H(Y_1,\dots,Y_{n-1})f(X_n)\\
&\qquad = EH(Y_1,\dots,Y_{n-1})f(X_n)\int_{\mathbb{S}} g(u)\,dP(X_0\in du|X_n,Y_{[1,n-1]})\\
&\qquad = EH(Y_1,\dots,Y_{n-1})\int_{\mathbb{S}}\int_{\mathbb{S}}\int_{\mathbb{S}} f(x)g(u)\,dP(X_0\in du|X_n=x,Y_{[1,n-1]})\\
&\qquad\qquad\times \lambda(v,x)\pi_{n-1}(v)\psi(dx)\psi(dv).
\end{aligned}\tag{3.3}$$

By arbitrariness of $g,f,H$ we derive from (3.2) and (3.3) that $(\psi\times\psi)$-a.s.

$$\begin{aligned}
&dP(X_0\in du|X_n,Y_{[1,n]})\int_{\mathbb{S}}\lambda(v;x)\pi_{n-1}(v)\psi(dv)\\
&\qquad = \psi(du)\int_{\mathbb{S}}\lambda(x',x)\rho_{n-1}(u;x')\pi_{n-1}(x')\psi(dx')
\end{aligned}\tag{3.4}$$

and the proof is complete. □

3.2. *Upper bound for* $|\rho_n(u,v')-\rho_n(u,v'')|$. Set

$$\overline{\rho}_n(u)=\operatorname*{ess\,sup}_{x\in\mathbb{S}}\rho_n(u,x)\quad\text{and}\quad \underline{\rho}_n(u)=\operatorname*{ess\,inf}_{x\in\mathbb{S}}\rho_n(u,x)\tag{3.5}$$

and introduce $\delta_n(u)=\overline{\rho}_n(u)-\underline{\rho}_n(u)$.

Lemma 3.2. *Under the assumptions of Theorem 1.1,*

$$\delta_1(u)\le\frac{(\lambda^*)^2}{\vartheta_*\lambda_\diamond}\vartheta(u),$$

$$\delta_n(u)\le\frac{(\lambda^*)^2}{\vartheta_*\lambda_\diamond}\vartheta(u)\exp\left\{-\frac{1}{\lambda^*}\sum_{k=2}^n\int\pi_{k-1}(x')\operatorname*{ess\,inf}_r\lambda(x',r)\psi(dx')\right\},\qquad n\ge 2.$$

Proof. Obviously,

$$\rho_1(u,x)\le\frac{(\lambda^*)^2\vartheta(u)}{\vartheta_*\int_{\mathbb{S}}\lambda(v,x)m(v)\psi(dv)}\le\frac{(\lambda^*)^2}{\vartheta_*\lambda_\diamond}\vartheta(u)$$

and, whereas $\delta_1(u) \le \rho_1(u,x)$, the first statement is valid. Further, by (3.1) we have

$$\begin{aligned}
&\rho_n(u,v') - \rho_n(u,v'') \\
&\quad = \frac{\int \pi_{n-1}(x')\lambda(x',v')\rho_{n-1}(u,x')\psi(dx')}{\int \pi_{n-1}(v)\lambda(v,v')\psi(dv)} \\
&\qquad - \frac{\int \pi_{n-1}(x')\lambda(x',v'')\rho_{n-1}(u,x')\psi(dx')}{\int \pi_{n-1}(v)\lambda(v,v'')\psi(dv)}.
\end{aligned} \tag{3.6}$$

Using the identities

$$\overline{\rho}_{n-1}(u) = \frac{\int \pi_{n-1}(x')\lambda(x',v')\overline{\rho}_{n-1}(u)\psi(dx')}{\int \pi_{n-1}(v)\lambda(v,v')\psi(dv)},$$

$$\underline{\rho}_{n-1}(u) = \frac{\int \pi_{n-1}(x')\lambda(x',v')\underline{\rho}_{n-1}(u)\psi(dx')}{\int \pi_{n-1}(v)\lambda(v,v')\psi(dv)},$$

rewrite (3.6) as

$$\begin{aligned}
&\rho_n(u,v') - \rho_n(u,v'') \\
&\quad = \overline{\rho}_{n-1}(u) - \underline{\rho}_{n-1}(u) \\
&\qquad - \frac{\int \pi_{n-1}(x')\lambda(x',v')[\overline{\rho}_{n-1}(u) - \rho_{n-1}(u,x')]\psi(dx')}{\int \pi_{n-1}(v)\lambda(v,v')\psi(dv)} \\
&\qquad - \frac{\int \pi_{n-1}(x')\lambda(x',v'')[\rho_{n-1}(u,x') - \underline{\rho}_{n-1}(u)]\psi(dx')}{\int \pi_{n-1}(v)\lambda(v,v'')\psi(dv)}.
\end{aligned} \tag{3.7}$$

Now, following the proof of Theorem 2.1, introduce

$$\alpha_n(u,x') := \frac{\rho_n(u,x') - \underline{\rho}_n(u)}{\delta_n(u)} \quad \text{and} \quad 1 - \alpha_n(u,x') := \frac{\overline{\rho}_n(u) - \rho_n(u,x')}{\delta_n(u)}$$

and rewrite (3.7) into

$$\begin{aligned}
&\rho_n(u,v') - \rho_n(u,v'') \\
&\quad = \delta_{n-1}(u)\Big(1 - \frac{\int \pi_{n-1}(x')\lambda(x',v')[1 - \alpha_{n-1}(u,x')]\psi(dx')}{\int \pi_{n-1}(v)\lambda(v,v')\psi(dv)} \\
&\qquad\qquad - \frac{\int \pi_{n-1}(x')\lambda(x',v'')\alpha_{n-1}(u,x')\psi(dx')}{\int \pi_{n-1}(v)\lambda(v,v'')\psi(dv)}\Big).
\end{aligned}$$

Since $0 \le \alpha_n(u,x') \le 1$, we have

$$\frac{\int \pi_{n-1}(x')\lambda(x',v')[1 - \alpha_{n-1}(u,x')]\psi(dx')}{\int \pi_{n-1}(v)\lambda(v,v')\psi(dv)}$$

$$
\begin{aligned}
&+\frac{\int \pi_{n-1}(x')\lambda(x',v'')\alpha_{n-1}(u,x')\psi(dx')}{\int \pi_{n-1}(v)\lambda(v,v'')\psi(dv)}\\
&\ge \int \pi_{n-1}(x')\frac{\lambda(x',v')}{\lambda^*}[1-\alpha_{n-1}(u,x')]\psi(dx')\\
&\quad+\int \pi_{n-1}(x')\frac{\lambda(x',v'')}{\lambda^*}\alpha_{n-1}(u,x')\psi(dx')\\
&\ge \frac{1}{\lambda^*}\int \pi_{n-1}(x')[\lambda(x',v')\wedge\lambda(x',v'')]\psi(dx')\\
&\ge \frac{1}{\lambda^*}\int \pi_{n-1}(x')\operatorname*{ess\,inf}_{r\in\mathbb{S}}\lambda(x',r)\psi(dx'),
\end{aligned}
$$

where the last inequality is valid $(\psi\times\psi)$-a.s. Hence,

$$
\begin{aligned}
&\rho_n(u,v')-\rho_n(u,v'')\\
&\qquad\le \delta_{n-1}(u)\Big(1-\frac{1}{\lambda^*}\int \pi_{n-1}(x')\operatorname*{ess\,inf}_{r}\lambda(x',r)\psi(dx')\Big)
\end{aligned}
$$

and, by arbitrariness of $v'$ and $v''$,

$$
\delta_n(u)\le \delta_{n-1}(u)\Big(1-\frac{1}{\lambda^*}\int \pi_{n-1}(x')\operatorname*{ess\,inf}_{r}\lambda(x',r)\psi(dx')\Big). \tag{3.8}
$$

The iteration of (3.8) and $\log(1-x)\le -x$, $0\le x<1$ provides the second statement. □

## 4. Proof of Theorem 1.1.

4.1. *Notation.* Let $(X^\nu,Y^\nu)$ and $(X^\beta,Y^\beta)$ be copies of $(X,Y)$ when the distribution density of $X_0$ is $\nu$ or $\beta$ and $Q^\nu$ and $Q^\beta$ be probability distributions of $(X^\nu,Y^\nu)$ and $(X^\beta,Y^\beta)$. Also let $\overline{Q}^\nu$ and $\overline{Q}^\beta$ be the distributions of $Y^\nu$ and $Y^\beta$ and $\overline{Q}^\nu_n$ and $\overline{Q}^\beta_n$ be the restrictions of $\overline{Q}^\nu$ and $\overline{Q}^\beta$ on the $\sigma$-algebras $Y^\nu_{[1,n]}$ and $Y^\beta_{[1,n]}$. Without loss of generality, we may assume that both copies are defined on the same probability space.

Set $\pi^\beta_n\langle f\rangle = E(f(X^\beta_n)|Y^\beta_{[1,n]})$ and let $\pi^\beta_n(x)$ be the corresponding filtering density.

4.2. *Absolute continuity $Q^\nu \ll Q^\beta$.* Since both Markov chains $(X^\nu,Y^\nu)$ and $(X^\beta,Y^\beta)$ have the same transition probability kernel, by the assumption $N\ll B$ we have

$$
Q^\nu\ll Q^\beta \qquad \text{with } \frac{dQ^\nu}{dQ^\beta}(X^\beta,Y^\beta)=\frac{\nu}{\beta}(X^\beta_0),
$$

$$\overline{Q}^\nu \ll \overline{Q}^\beta \qquad \text{with } \frac{d\overline{Q}^\nu}{d\overline{Q}^\beta}(Y^\beta) = E\Big(\frac{\nu}{\beta}(X_0^\beta)\Big|Y^\beta_{[1,\infty)}\Big),\ Y^\beta_{[1,\infty)} = \bigvee_{n\ge 1} Y^\beta_{[1,n]},$$

$$\overline{Q}^\nu_n \ll \overline{Q}^\beta_n \qquad \text{with } \frac{d\overline{Q}^\nu_n}{d\overline{Q}^\beta_n}(Y^\beta) = E\Big(\frac{\nu}{\beta}(X_0^\beta)\Big|Y^\beta_{[1,n]}\Big),\ n \ge 1.$$

4.3. *Proof of* (1.4). It suffices to show that

$$E\Big(\frac{\nu}{\beta}(X_0^\beta)\Big|\bigcap_{n\ge 0}\{Y^\beta_{[1,\infty)} \vee X^\beta_{[n,\infty)}\}\Big)$$

is a $Y^\beta_{[1,\infty)}$-measurable random variable. Indeed, then we have a.s.

$$\begin{aligned}
&E\Big(\frac{\nu}{\beta}(X_0^\beta)\Big|\bigcap_{n\ge 0}\{Y^\beta_{[1,\infty)} \vee X^\beta_{[n,\infty)}\}\Big)\\
&\quad = E\Big(E\Big(\frac{\nu}{\beta}(X_0^\beta)\Big|\bigcap_{n\ge 0}\{Y^\beta_{[1,\infty)} \vee X^\beta_{[n,\infty)}\}\Big)\Big|Y^\beta_{[1,\infty)}\Big)\\
&\quad = E\Big(\frac{\nu}{\beta}(X_0^\beta)\Big|Y^\beta_{[1,\infty)}\Big).
\end{aligned}$$

In order to verify the required $Y^\beta_{[1,\infty)}$-measurability, notice first that by the reverse martingale convergence theorem,

$$E\Big(\frac{\nu}{\beta}(X_0^\beta)\Big|\bigcap_{n\ge 0}\{Y^\beta_{[1,\infty)} \vee X^\beta_{[n,\infty)}\}\Big) = \lim_{n\to\infty} E\Big(\frac{\nu}{\beta}(X_0^\beta)\Big|Y^\beta_{[1,\infty)} \vee X^\beta_{[n,\infty)}\Big)$$

and by the Markov property,

$$\begin{aligned}
E\Big(\frac{\nu}{\beta}(X_0^\beta)\Big|Y^\beta_{[1,\infty)} \vee X^\beta_{[n,\infty)}\Big) &= E\Big(\frac{\nu}{\beta}(X_0^\beta)\Big|Y^\beta_{[1,n]} \vee X^\beta_n\Big)\\
&\equiv \int_{\mathbb{S}} \frac{\nu}{\beta}(u)\rho_n(u, X_n^\beta)\psi(du).
\end{aligned}$$

So by the aforementioned martingale convergence theorem,

$$\lim_{n\to\infty}\int_{\mathbb{S}} \frac{\nu}{\beta}(u)\rho_n(u, X_n^\beta)\psi(du) \qquad \text{exists a.s.}$$

On the other hand, since with $\overline{\rho}_n(u)$ and $\underline{\rho}_n(u)$ defined in (3.5) we have

$$\int_{\mathbb{S}} \frac{\nu}{\beta}(u)\underline{\rho}_n(u)\psi(du) \le \int_{\mathbb{S}} \frac{\nu}{\beta}(u)\rho_n(u, X_n^\beta)\psi(du) \le \int_{\mathbb{S}} \frac{\nu}{\beta}(u)\overline{\rho}_n(u)\psi(du)$$

and, by Lemma 3.2 with $\vartheta(x) \equiv \beta(x)$, for $n \ge 2$,

$$\int_{\mathbb{S}} \frac{\nu}{\beta}(u)[\overline{\rho}_n(u) - \underline{\rho}_n(u)]\psi(du) \le \frac{(\lambda^*)^2}{\beta_* \lambda_\diamond} \exp\left\{-\frac{1}{\lambda^*}\sum_{k=2}^{n} \int_{\mathbb{S}} \pi_{k-1}^{\beta}(x') \operatorname*{ess\,inf}_r \lambda(x',r)\psi(dx')\right\}.$$

By Theorem 2.2, $\lim_{n\to\infty} n^{-1}\sum_{k=2}^{n}\int \pi_{k-1}^{\beta}(x')\operatorname{ess\,inf}_r \lambda(x',r)\psi(dx') = \lambda_\diamond\,(>0)$, a.s. Consequently, $\lim_{n\to\infty}\int_{\mathbb{S}}[\overline{\rho}_n(u) - \underline{\rho}_n(u)]\nu(u)/\beta(u)\psi(du) = 0$, a.s. and, therefore,

$$E\left(\frac{\nu}{\beta}(X_0^{\beta})\middle| \bigcap_{n\ge 0}\{Y_{[1,\infty)}^{\beta} \vee X_{[n,\infty)}^{\beta}\}\right) = \lim_{n\to\infty}\int_{\mathbb{S}} \frac{\nu}{\beta}(u)\overline{\rho}_n(u)\psi(du).$$

The desired statement holds since $\int_{\mathbb{S}} \overline{\rho}_n(u)\nu(u)/\beta(u)\psi(du)$ is a $Y_{[1,\infty)}$-measurable random variable for any $n$.

4.4. *Proof of* (1.8). For fixed $n$, let $\phi^{\nu}(y_{[1,n]})$ and $\phi^{\beta}(y_{[1,n]})$ be measurable functions of the arguments $y_1,\dots,y_n$, $y_i \in \mathbb{R}^p$, such that for fixed bounded and measurable function $f$,

$$\phi^{\nu}(Y_{[1,n]}^{\nu}) = \pi_n^{\nu}\langle f\rangle \quad \text{and} \quad \phi^{\beta}(Y_{[1,n]}^{\beta}) = \pi_n^{\beta}\langle f\rangle.$$

Recall that $N \ll B$ provides $Q^{\nu} \ll Q^{\beta}$. So we may define $\pi_n^{\beta\nu}\langle f\rangle = \phi^{\beta}(Y_{[1,n]}^{\nu})$. Since $Q^{\beta} \ll Q^{\nu}$ is not assumed, we set

$$\pi_n^{\nu\beta}\langle f\rangle = \phi^{\nu}(Y_{[1,n]}^{\beta}) I\left(\frac{d\overline{Q}_n^{\nu}}{d\overline{Q}_n^{\beta}}(Y_{[1,n]}^{\beta}) > 0\right).$$

Lemma 4.1. *Under the assumptions of Theorem 1.1, for $n \ge 2$,*

$$\frac{d\overline{Q}_n^{\nu}}{d\overline{Q}_n^{\beta}}(Y^{\beta})\|\pi_n^{\nu\beta} - \pi_n^{\beta}\| \le \frac{(\lambda^*)^2}{\beta_*\lambda_\diamond}\exp\left\{-\frac{1}{\lambda^*}\sum_{k=2}^{n}\int \pi_{k-1}^{\beta}(x')\operatorname*{ess\,inf}_r \lambda(x',r)\psi(dx')\right\}. \tag{4.1}$$

Proof. We show first that

$$\frac{d\overline{Q}_n^{\nu}}{d\overline{Q}_n^{\beta}}(Y^{\beta})(\pi_n^{\nu\beta}(x) - \pi_n^{\beta}(x)) = \pi_n^{\beta}(x)\left[\int_{\mathbb{S}}\frac{\nu}{\beta}(u)\rho_n(u,x)\psi(du) - \int_{\mathbb{S}}\int_{\mathbb{S}}\frac{\nu}{\beta}(u)\rho_n(u,x')\pi_n^{\beta}(x')\psi(dx')\psi(du)\right]. \tag{4.2}$$

Let $g$ and $f$ be bounded measurable functions of arguments $y_1,\dots,y_n$, $y_i\in\mathbb{R}^d$ and $x\in\mathbb{S}$. Then,

$$\begin{aligned}
&E\Big[E\Big(\frac{\nu}{\beta}(X_0^\beta)\Big|Y^\beta_{[1,n]}\Big)\pi_n^{\nu\beta}\langle f\rangle - E\Big(\frac{\nu}{\beta}(X_0^\beta)f(X_n^\beta)\Big|Y^\beta_{[1,n]}\Big)\Big]g(Y^\beta_{[1,n]})\\
&\qquad = E\frac{\nu}{\beta}(X_0^\beta)g(Y^\beta_{[1,n]})[\pi_n^{\nu\beta}\langle f\rangle - f(X_n^\beta)]\\
&\qquad = Eg(Y^\nu_{[1,n]})[\pi_n^\nu\langle f\rangle - f(X_n^\nu)]\\
&\qquad = 0
\end{aligned}$$

provides

$$\frac{d\overline{Q}{}^\nu_n}{d\overline{Q}{}^\beta_n}(Y^\beta)\pi_n^{\nu\beta}\langle f\rangle = E\Big(\frac{\nu}{\beta}(X_0^\beta)f(X_n^\beta)\Big|Y^\beta_{[1,n]}\Big). \tag{4.3}$$

Further, by (4.3),

$$\begin{aligned}
&\frac{d\overline{Q}{}^\nu_n}{d\overline{Q}{}^\beta_n}(Y^\beta)(\pi_n^{\nu\beta}\langle f\rangle - \pi_n^\beta\langle f\rangle)\\
&\quad = E\Big(\frac{\nu}{\beta}(X_0^\beta)f(X_n^\beta)\Big|Y^\beta_{[1,n]}\Big) - E\Big(\frac{\nu}{\beta}(X_0^\beta)\Big|Y^\beta_{[1,n]}\Big)E(f(X_n^\beta)|Y^\beta_{[1,n]})\\
&\quad = E\Big(f(X_n^\beta)\Big[E\Big\{\frac{\nu}{\beta}(X_0^\beta)\Big|X_n^\beta,Y^\beta_{[1,n]}\Big\} - E\Big\{\frac{\nu}{\beta}(X_0^\beta)\Big|Y^\beta_{[1,n]}\Big\}\Big]\Big|Y^\beta_{[1,n]}\Big)
\end{aligned}$$

and so, due to arbitrariness of $f$, (4.2) holds true.

The right-hand side of (4.2) is evaluated from above as follows:

$$\begin{aligned}
&\Big|\pi_n^\beta(x)\Big[\int_{\mathbb{S}}\frac{\nu}{\beta}(u)\rho_n(u,x)\psi(du)\\
&\qquad\qquad - \int_{\mathbb{S}}\int_{\mathbb{S}}\frac{\nu}{\beta}(u)\rho_n(u,x')\pi_n^\beta(x')\psi(dx')\psi(du)\Big]\Big|\\
&\quad = \Big|\pi_n^\beta(x)\Big[\int_{\mathbb{S}}\int_{\mathbb{S}}\frac{\nu}{\beta}(u)[\rho_n(u,x)-\rho_n(u,x')]\pi_n^\beta(x')\psi(du)\psi(dx')\Big]\Big|\\
&\quad \le \pi_n^\beta(x)\int_{\mathbb{S}}\frac{\nu}{\beta}(u)\delta_n(u)\psi(du)
\end{aligned}$$

and the desired result follows by Lemma 3.2 with $\vartheta=\beta$. $\square$

By (4.1),

$$\begin{aligned}
&\frac{1}{n}\log\frac{d\overline{Q}{}^\nu_n}{d\overline{Q}{}^\beta_n}(Y^\beta)\|\pi_n^{\nu\beta}-\pi_n^\beta\|\\
&\qquad \le \frac{1}{n}\log\frac{(\lambda^*)^2}{\beta_*\lambda_\diamond} - \frac{1}{\lambda^* n}\sum_{k=2}^n\int_{\mathbb{S}}\pi_{k-1}^\beta(x')\operatorname*{ess\,inf}_{r\in\mathbb{S}}\lambda(x',r)\psi(dx')
\end{aligned}$$

and, since by Theorem 2.2,

$$\lim_{n\to\infty}\frac{1}{n}\sum_{k=2}^{n}\int_{\mathbb{S}}\pi^{\beta}_{k-1}(x')\operatorname*{ess\,inf}_{r}\lambda(x',r)\psi(dx')$$
$$=\int_{\mathbb{S}}m(x')\operatorname*{ess\,inf}_{r}\lambda(x',r)\psi(dx')=\lambda_{\diamond},$$

we get

$$P\Big(\limsup_{n\to\infty}\frac{1}{n}\log\Big\{\frac{d\overline{Q}^{\nu}_{n}}{d\overline{Q}^{\beta}_{n}}(Y^{\beta})\|\pi^{\nu\beta}_{n}-\pi^{\beta}_{n}\|\Big\}\le-\frac{\lambda_{\diamond}}{\lambda^{*}}\Big)=1.$$

On the other hand, owing to $\overline{Q}^{\nu}\ll\overline{Q}^{\beta}$,

$$P\Big(\limsup_{n}\frac{1}{n}\log\frac{d\overline{Q}^{\nu}_{n}}{d\overline{Q}^{\beta}_{n}}(Y^{\nu})\|\pi^{\nu}_{n}-\pi^{\beta\nu}_{n}\|\le-\frac{\lambda_{\diamond}}{\lambda^{*}}\Big)=1.$$

Then, the statement of Theorem 1.1 is valid, provided that

$$\lim_{n\to\infty}\frac{1}{n}\log\frac{d\overline{Q}^{\nu}_{n}}{d\overline{Q}^{\beta}_{n}}(Y^{\nu})=0. \tag{4.4}$$

Notice that

$$\frac{d\overline{Q}^{\nu}_{n}}{d\overline{Q}^{\beta}_{n}}(Y^{\beta})=E\Big(\frac{d\overline{Q}^{\nu}}{d\overline{Q}^{\beta}}(Y^{\beta})\Big|Y^{\beta}_{[1,n]}\Big),$$

that is, $\frac{d\overline{Q}^{\nu}_{n}}{d\overline{Q}^{\beta}_{n}}(Y^{\beta})$ is a nonnegative uniformly integrable martingale with the limit point $\frac{d\overline{Q}^{\nu}}{d\overline{Q}^{\beta}}(Y^{\beta})$ such that $0\le\frac{d\overline{Q}^{\nu}}{d\overline{Q}^{\beta}}(Y^{\beta})<\infty$, a.s. Again taking into account $\overline{Q}^{\nu}\ll\overline{Q}^{\beta}$, we conclude that

$$P\Big(\lim_{n\to\infty}\frac{d\overline{Q}^{\nu}_{n}}{d\overline{Q}^{\beta}_{n}}(Y^{\nu})=\frac{d\overline{Q}^{\nu}}{d\overline{Q}^{\beta}}(Y^{\nu})\Big)=1$$

with

$$P\Big(0\le\frac{d\overline{Q}^{\nu}_{n}}{d\overline{Q}^{\beta}_{n}}(Y^{\beta})<\infty\Big)=1.$$

Moreover, on account of

$$P\Big(\frac{d\overline{Q}^{\nu}}{d\overline{Q}^{\beta}}(Y^{\nu})=0\Big)=EI\Big(\frac{d\overline{Q}^{\nu}}{d\overline{Q}^{\beta}}(Y^{\beta})=0\Big)\frac{d\overline{Q}^{\nu}}{d\overline{Q}^{\beta}}(Y^{\beta})=0,$$

we have that $P(0<\frac{d\overline{Q}^{\nu}}{d\overline{Q}^{\beta}}(Y^{\nu})<\infty)=1$.

Thus, (4.4) is valid.

**5. Discussion and one counterexample.** In this paper the exponential stability of the filtering equation is verified for a class of ergodic Markov signals whose transition probability density may vanish on a part of the state space.

To the best of our knowledge, the assumption of Theorem 1.1 is the weakest currently known one, when the filter is stable regardless of the observations.

The obtained result, as well as the result based on the mixing conditions, is ideally compatible with compact space $\mathbb{S}$. Except for several known examples (see, e.g., Example 1 in [6]) with $\mathbb{S}=\mathbb{R}$, these conditions are not applicable in many practical filtering schemes with noncompact $\mathbb{S}$. Particularly, the Kalman–Bucy filter is out of their scope, though its stability is well known (see [3, 7, 12, 13, 15] and [11], Sections 14.6 and 16.2).

Since both types of conditions imply the geometric ergodicity of the signal, it may seem that it provides the filtering stability. Example 5.1 demonstrates that this is false in general.

EXAMPLE 5.1. Following [8], we consider Markov chain $X_n^\nu$, with values in $\mathbb{S}=\{1,2,3,4\}$ and the transition probability matrix

$$\Lambda=\begin{pmatrix}0.5 & 0.5 & 0 & 0\\ 0 & 0.5 & 0.5 & 0\\ 0 & 0 & 0.5 & 0.5\\ 0.5 & 0 & 0 & 0.5\end{pmatrix}.$$

Owing to $\lambda_*=0$ and $\lambda_\diamond=0$, neither mixing nor our conditions are not applicable. On the other hand, since 1 is the simple eigenvalue, the geometric ergodicity holds for the unique invariant distribution $(0.25;0.25;0.25;0.25)$.

Suppose that

$$Y_n^\nu=I(X_n^\nu=1)+I(X_n^\nu=3),\qquad n\ge 1.$$

In the context of our setting, the measure $\psi$ is defined in an obvious way (see Example 1.1) and $\nu=(\nu_1;\nu_2;\nu_3;\nu_4)$ and $\beta=(\beta_1;\beta_2;\beta_3;\beta_4)$ are corresponding densities with respect to $\psi$.

In [8], Kaijser notes that the vector process $\pi_n$ with the entries $\pi_n(i)=P(X_n=i|Y_{[1,n]})$, $i=1,2,3,4$ is not ergodic. We show also that

$$\|\pi_n^{\beta\nu}-\pi_n^\nu\|=\sum_{k=1}^4|\pi_n^{\beta\nu}(k)-\pi_n^\nu(k)|$$

remains strictly positive for any $n\ge 1$.

In this case, the recursion (1.2) reads $(n\ge 1)$,

$$\pi_n^\nu(1)=[\pi_{n-1}^\nu(1)+\pi_{n-1}^\nu(4)]Y_n^\nu,$$

$$\pi_n^\nu(2) = [\pi_{n-1}^\nu(2) + \pi_{n-1}^\nu(1)](1 - Y_n^\nu),$$
$$\pi_n^\nu(3) = [\pi_{n-1}^\nu(3) + \pi_{n-1}^\nu(2)]Y_n^\nu,$$
$$\pi_n^\nu(4) = [\pi_{n-1}^\nu(4) + \pi_{n-1}^\nu(3)](1 - Y_n^\nu)$$

subject to $\pi_0^\nu(k) = \nu_k$, $k = 1, \dots, 4$. Similarly,

$$\pi_n^{\beta\nu}(1) = [\pi_{n-1}^{\beta\nu}(1) + \pi_{n-1}^{\beta\nu}(4)]Y_n^\nu,$$
$$\pi_n^{\beta\nu}(2) = [\pi_{n-1}^{\beta\nu}(2) + \pi_{n-1}^{\beta\nu}(1)](1 - Y_n^\nu),$$
$$\pi_n^{\beta\nu}(3) = [\pi_{n-1}^{\beta\nu}(3) + \pi_{n-1}^{\beta\nu}(2)]Y_n^\nu,$$
$$\pi_n^{\beta\nu}(4) = [\pi_{n-1}^{\beta\nu}(4) + \pi_{n-1}^{\beta\nu}(3)](1 - Y_n^\nu)$$

subject to $\pi_0^{\beta\nu}(k) = \beta_k$, $k = 1, \dots, 4$.

By virtue of the obvious identities,

$$\pi_n^\nu(k)Y_n^\nu = \pi_n^\nu(k) \quad \text{and} \quad \pi_n^{\beta\nu}(k)Y_n^\nu = \pi_n^{\beta\nu}(k), \qquad k = 1, 3,$$
$$\pi_n^\nu(k)(1 - Y_n^\nu) = \pi_n^\nu(k) \quad \text{and} \quad \pi_n^{\beta\nu}(k)(1 - Y_n^\nu) = \pi_n^{\beta\nu}(k), \qquad k = 2, 4,$$

we find

$$\begin{aligned}|\pi_n^\nu(1) - \pi_n^{\beta\nu}(1)| &= |\pi_{n-1}^\nu(1) - \pi_{n-1}^{\beta\nu}(1) + \pi_{n-1}^\nu(4) - \pi_{n-1}^{\beta\nu}(4)|Y_n^\nu \\ &= |\pi_{n-1}^\nu(1) - \pi_{n-1}^{\beta\nu}(1)|Y_{n-1}^\nu Y_n^\nu \\ &\quad + |\pi_{n-1}^\nu(4) - \pi_{n-1}^{\beta\nu}(4)|(1 - Y_{n-1}^\nu)Y_n^\nu.\end{aligned}$$

Similarly,

$$\begin{aligned}|\pi_n^\nu(2) - \pi_n^{\beta\nu}(2)| &= |\pi_{n-1}^\nu(2) - \pi_{n-1}^{\beta\nu}(2)|(1 - Y_{n-1}^\nu)(1 - Y_n^\nu) \\ &\quad + |\pi_{n-1}^\nu(1) - \pi_{n-1}^{\beta\nu}(1)]Y_{n-1}^\nu(1 - Y_n^\nu), \\ |\pi_n^\nu(3) - \pi_n^{\beta\nu}(3)| &= |\pi_{n-1}^\nu(3) - \pi_{n-1}^{\beta\nu}(3)|Y_{n-1}^\nu Y_n^\nu \\ &\quad + |\pi_{n-1}^\nu(2) - \pi_{n-1}^{\beta\nu}(2)|(1 - Y_{n-1}^\nu)Y_n^\nu, \\ |\pi_n^\nu(4) - \pi_n^{\beta\nu}(4)| &= |\pi_{n-1}^\nu(4) - \pi_{n-1}^{\beta\nu}(4)|(1 - Y_{n-1}^\nu)(1 - Y_n^\nu) \\ &\quad + |\pi_{n-1}^\nu(3) - \pi_{n-1}^{\beta\nu}(3)]Y_{n-1}^\nu(1 - Y_n^\nu).\end{aligned} \tag{5.1}$$

Now, (5.1) provides $\|\pi_n^\nu - \pi_n^{\beta\nu}\| \equiv \|\pi_{n-1}^\nu - \pi_{n-1}^{\beta\nu}\|$, $n \ge 2$.

Hence, $\|\pi_n^\nu - \pi_n^{\beta\nu}\| \equiv \|\pi_1^\nu - \pi_1^{\beta\nu}\|$. Set

$$c_1 = |\nu_1 - \beta_1 + \nu_4 - \beta_4| + |\nu_3 - \beta_3 + \nu_2 - \beta_2|,$$
$$c_2 = |\nu_2 - \beta_2 + \nu_1 - \beta_1| + |\nu_4 - \beta_4 + \nu_3 - \beta_3|$$

and assume that $c_1 \wedge c_2 > 0$. Since

$$\|\pi_1^\nu - \pi_1^{\beta\nu}\| = c_1 Y_1^\nu + c_2(1 - Y_1^\nu) \geq c_1 \wedge c_2,$$

we have $\|\pi_n^\nu - \pi_n^{\beta\nu}\| \geq c_1 \wedge c_2$.

**Acknowledgements.** The authors gratefully acknowledge the anonymous referees whose comments and advice allowed us to improve the paper significantly.

Department of Electrical
Engineering-Systems
Tel Aviv University
69978 Tel Aviv
Israel
e-mail: pavelm@eng.tau.ac.il
e-mail: liptser@eng.tau.ac.il